\documentclass[a4paper, 10pt]{article}
 
\usepackage{enumerate}
\usepackage{parselines} 
\usepackage{authblk}
\usepackage[latin1]{inputenc}
\usepackage[francais]{babel}
\usepackage[T1]{fontenc}
\usepackage{amsfonts}
\usepackage{graphicx}
\usepackage{ulem}
\usepackage{lmodern}
\usepackage[top=3cm, bottom=3cm, left=4cm, right=4cm]{geometry}
\usepackage{amsthm}
\usepackage{amsmath}
\usepackage{amssymb}
\usepackage{mathrsfs}
\usepackage{enumerate}
\usepackage{sectsty}

\newtheorem{Lemme}{Lemme}
\newtheorem{Theoreme}{Théorème} 

\usepackage{fancyhdr}
\allsectionsfont{\centering}

\begin{document}

\title{SUR L'ANNULATION DE LA VALEUR CENTRALE DE LA FONCTION $L$ DE HECKE}
\author{Quentin Gazda\footnote{quentin.gazda@ens-paris-saclay.fr}}
\affil{\'Ecole Normale Sup\'erieure Paris-Saclay}

\maketitle

\textsl{
\begin{flushright}
\`A la mémoire de Wendy
\end{flushright}
}
\par
\par

\begin{abstract}

Dans cette note, en se basant sur des résultats de Henri Cohen et Winfried Kohnen, nous démontrons que pour tout entier $k$ supérieur ou égal à $12$ et divisible par $4$, il existe une forme parabolique propre de poids $k$ pour $\operatorname{SL}_2(\textbf{Z})$ telle que sa fonction $L$ de Hecke ne s'annule pas en $k/2$. \\

In this note, following results from Henri Cohen and Winfried Kohnen, we show that for all integer $k$ greater than $12$ and divisible by $4$, there exists a cuspidal eigenform of weight $k$ for the full modular group $\operatorname{SL}_2(\textbf{Z})$ such that its Hecke $L$-function does not vanish on $k/2$.
\end{abstract}

\section{Introduction et énoncé du résultat}
Soit $f$ une forme parabolique de poids $k$ entier (pair) pour le groupe modulaire $\Gamma(1):=\operatorname{SL}_2(\textbf{Z})$, et soit $L(f,s)$ ($\operatorname{Re}(s)>k/2+1$) la fonction $L$ de Hecke associée. Il est connu, depuis Hecke, que cette fonction $L$ se prolonge holomorphiquement à l'ensemble du plan complexe en une fonction satisfaisant à l'équation fonctionnelle :
\begin{equation} \label{eqfonc}
(2\pi)^{-s}\Gamma(s)L(f,s)=(-1)^{k/2}(2\pi)^{-(k-s)}\Gamma(k-s)L(f,k-s) \quad (s\in \textbf{C}).
\end{equation}
Il résulte de (\ref{eqfonc}) que $L(f,k/2)=(-1)^{k/2}L(f,k/2)$ (prendre $s=k/2$); ainsi, lorsque $k\equiv 2 \pmod{4}$, l'équation fonctionnelle (\ref{eqfonc}) force $L(f,k/2)$ à être nulle. Seulement, rien ne semble indiquer un résultat similaire dans le cas où $k\equiv 0 \pmod{4}$. Nous allons effectivement montrer que lorsque $4$ divise $k$ et que l'espace des formes cuspidales de poids $k$ pour $\Gamma(1)$, noté $\mathcal{S}_{k}(\Gamma(1))$, est non réduit à l'espace nul, il existe une forme parabolique $f$ telle que $L(f,k/2)$ soit non nulle.
\begin{Theoreme} \label{tun}
Soit $k$ un entier. Il existe une forme propre de Hecke de poids $k$ (pour $\Gamma(1)$) telle que $L(f,k/2) \neq 0$ si, et seulement si $k$ est supérieur ou égal à $12$ et divisible par $4$. 
\end{Theoreme}

Afin de démontrer ce résutat, nous allons étudier la forme parabolique duale de la forme linéaire sur $\mathcal{S}_k(\Gamma(1))$ qui à $\varphi$ associe le nombre $L(\varphi,k/2)$ (vis-à-vis du produit scalaire de Petersson). Cette forme fut introduite dans un cadre plus général par Cohen (dans \cite{Cohen}) pour trouver des bases de l'espace des formes paraboliques, puis par Kohnen (dans \cite{Kohnen}) afin d'obtenir une version faible de l'hypothèse de Riemann généralisée. La non nullité de son premier coefficient de Fourier nous indiquera la non nullité de la forme linéaire et donc l'existence de la forme cuspidale propre recherchée.
\section{Preuve du Théorème \ref{tun}}
Commençons par montrer que lorsque l'entier $k$ n'est pas de la forme indiquée au Théorème \ref{tun} alors $L(f,k/2)$ est nul. L'espace $\mathcal{S}_k(\Gamma(1))$ étant nul pour $k$ impair ou $k<12$, il suffit de se contenter des cas où $k\equiv 2 \pmod 4$ et $k\geq 12$. Mais alors (et c'est essentiellement la remarque faite en introduction), la fonction $L$ de Hecke de $f\in \mathcal{S}_k(\Gamma(1))$ satisfait à l'équation fonctionnelle (\ref{eqfonc}) et ainsi $L(f,k/2)=0$. \\

Soient $k$ un entier supérieur ou égal à $12$ divisible par $4$, $\operatorname{H}_k$ l'ensemble des formes propres de Hecke de poids $k$ pour $\Gamma(1)$, $(\cdot,\cdot)_k$ le produit scalaire de Petersson de poids $k$ et $\|\cdot \|_k$ la norme associée. Il est connu que $\operatorname{H}_k$ est une base orthonorgonale de $\mathcal{S}_k(\Gamma(1))$. Nous allons montrer que la fonction $L$ de Hecke associée à la forme cuspidale suivante :
\begin{equation}
f_k:=\sum_{f\in \operatorname{H}_k}{\frac{f}{\|f\|_k^2}}, \nonumber
\end{equation}
ne s'annule pas en $k/2$, ce qui impliquera en particulier qu'il existe $f\in\operatorname{H}_k$ telle que sa fonction $L$ ne s'annule pas en $k/2$. Pour cela, définissons pour $\tau \in \mathcal{H}$ (un élément du demi-plan de Poincaré) :
\begin{equation}
R_k(\tau):=c_k\sum_{\gamma \in \Gamma(1)}{(a\tau+b)^{-\frac{k}{2}}(c\tau+d)^{-\frac{k}{2}}} \quad \left(\text{avec}~\gamma=\begin{pmatrix} a & b \\ c & d  \end{pmatrix}\right) \nonumber
\end{equation}
où $c_k:=(-1)^{\frac{k}{4}}(8\pi)^{\frac{k}{2}-1}\left(\frac{k}{2}-1\right)!/(k-2)!$. La fonction qui à $\tau$ associe $R_k(\tau)$ fut étudiée (au coefficient $c_k$ près) indépendamment par Cohen et Kohnen : elle définit notamment un élément de $\mathcal{S}_k(\Gamma(1))$ (rappelons que $k\geq 12$). On retrouve chez les deux auteurs le lemme suivant (voir par exemple le lemme $1$ dans \cite{Kohnen}) :
\begin{Lemme} 
Soit $f$ un élément de $\mathcal{S}_k(\Gamma(1))$. Alors, $(f,R_k)=L(f,k/2)$. En particulier, nous avons l'égalité :
\begin{equation} \label{dual}
R_k=\sum_{f\in \operatorname{H}_k}{\frac{f}{\|f\|_k^2}L\left(f,\frac{k}{2}\right)}.
\end{equation}
\end{Lemme}

Notons $r_k(n)$ les coefficients de Fourier de $R_k$ de sorte à ce que
\begin{equation}
R_k(\tau)=\sum_{n=1}^{\infty}{r_k(n)q^n} \quad (q=e^{2i\pi \tau}). \nonumber
\end{equation}
Nous obtenons alors $r_k(1)=L(f_k,k/2)$ en identifiant le premier coefficient de Fourier de chaque membre de l'égalité (\ref{dual}). Il suffit ainsi de montrer que $r_k(1)$ est non nul. Nous disposons d'un second lemme (voir le corollaire $3.2$ dans \cite{Cohen}) :
\begin{Lemme}
Pour $n$ un entier strictement positif et $k$ comme précédemment, nous avons
\begin{equation}
r_{k}(n)=\frac{(8\pi)^{\frac{k}{2}-1}}{4(k-2)!}n^{\frac{k}{2}-1}\left(1+(-1)^{k/4}\sqrt{2\pi}\sum_{m=1}^{\infty}{\gamma_n(m)\sqrt{\frac{n\pi}{m}}J_{\frac{k-1}{2}}\left(\frac{n\pi}{m}\right)}\right) \nonumber
\end{equation}
où, pour tout entiers positifs $m$ et $n$, $\gamma_n(m)$ est définit par la somme finie :
\begin{equation}
\gamma_n(m)=\sum_{\substack{ac=m\\(a,c)=1}}{\cos\left[\pi n\left(\frac{a'}{c}-\frac{c'}{a}\right)\right]} \nonumber
\end{equation}
(la somme porte sur les couples $(a,c)$ d'entiers positifs tels que $ac=m$ et $a$, $c$ premiers entre eux) avec $a'$ (resp. $c'$) l'inverse de $a \pmod c$ (resp. $c \pmod a$) et $J_{\frac{k-1}{2}}$ la fonction de Bessel d'ordre $(k-1)/2$.
\end{Lemme}
Il s'agit donc de montrer que l'expression (proportionnelle à $r_k(1)$)
\begin{equation}
1+(-1)^{k/4}\sqrt{2\pi}\sum_{m=1}^{\infty}{\gamma_1(m)\sqrt{\frac{\pi}{m}}J_{\frac{k-1}{2}}\left(\frac{\pi}{m}\right)} \nonumber
\end{equation}
est non nulle. Nous avons $|\gamma_1(m)|\leq d(m)$ (où $d(m)$ correspond au nombre de diviseurs positifs de $m$) ainsi que l'inégalité classique pour les fonctions de Bessel 
\begin{equation}
\left|J_{\frac{k-1}{2}}(x)\right|\leq \frac{1}{\Gamma\left(\frac{k+1}{2}\right)}\left(\frac{x}{2}\right)^{\frac{k-1}{2}}=\sqrt{\frac{2}{\pi}}\frac{(k/2)!}{k!}2^{\frac{k}{2}}x^{\frac{k-1}{2}} \quad (\text{pour}~x\geq 0)
\end{equation}
d'après les valeurs aux demi-entiers de la fonction Gamma. Il en résulte l'inégalité suivante :
\begin{equation}
\left|\sqrt{2\pi}\sum_{m=1}^{\infty}{\gamma_1(m)\sqrt{\frac{\pi}{m}}J_{\frac{k-1}{2}}\left(\frac{\pi}{m}\right)}\right| \leq 2(2\pi)^{\frac{k}{2}}\frac{(k/2)!}{k!}\left(\sum_{m=1}^{\infty}{\frac{d(m)}{m^{\frac{k}{2}}}}\right)=2(2\pi)^{\frac{k}{2}}\frac{(k/2)!}{k!}\zeta\left(\frac{k}{2}\right)^2. \nonumber
\end{equation}
Puisque $k\geq 12$, nous avons $\zeta(k/2)^2 \leq \zeta(6)^2$ ainsi que
\begin{equation}
(2\pi)^{\frac{k}{2}}\frac{(k/2)!}{k!}\leq \left(\frac{2\pi}{k/2+1}\right)\left(\frac{2\pi}{k/2+2}\right)^{\frac{k}{2}-1}\leq \left(\frac{2\pi}{7}\right)\left(\frac{2\pi}{8}\right)^5. \nonumber
\end{equation}
Finalement, 
\begin{equation}
\left|1+(-1)^{k/4}\sqrt{2\pi}\sum_{m=1}^{\infty}{\gamma_1(m)\sqrt{\frac{\pi}{m}}J_{\frac{k-1}{2}}\left(\frac{\pi}{m}\right)}\right|\geq 1-2\left(\frac{2\pi}{7}\right)\left(\frac{2\pi}{8}\right)^5\zeta(6)^2 >0, \nonumber
\end{equation}
et donc $|r_k(1)|>0$ ce qui conclut.

\section*{Remerciements}

Cette note fut écrite lors de mon stage de Master $1$ à l'\'ENS de Lyon encadré par François Brunault et Gabriel Dospinescu. Je les remercie vivement pour m'avoir guidé vers ce problème. Je souhaite également remercier l'ensemble de l'équipe de l'UMPA pour son accueille si sympathique et pour m'avoir proposé d'aussi bonnes conditions de travail.

\end{document}